# The Natural Philosophy of Kazuo Kondo

Grenville J. Croll
Alternative Natural Philosophy Association
grenvillecroll@gmail.com

*Kazuo Kondo (1911-2001) was Chair of the Department of Mathematical Engineering at the University of Tokyo, Japan. Over a period of 50 years, he and a few colleagues wrote and published a voluminous series of papers and monographs on the applications of analytical geometry within a diverse range of subjects in the natural sciences. Inspired by Otto Fischer's attempt at a quaternionic unified theory in the late 1950's he adopted the mathematics of the revered Akitsugu Kawaguchi to produce his own speculative unified theory. The theory appears to successfully apply Kawaguchi's mathematics to the full range of natural phenomena, from the structure of fundamental particles to the geometry of living beings. The theories are testable and falsifiable. Kondo and his theories are now almost completely unknown and this paper serves as the barest introduction to his work.* v2.39

## 1 Biographical

Very little was known of Kazuo Kondo. He was born on the 2nd January 1911 and died of lung cancer on the 1st December 2001. He was known to be married and his wife Reiko of 37 years survived him, but there is no documented evidence of any children or other family members.

After graduating from the Department of Aeronautical Engineering, Faculty of Engineering, Tokyo Imperial University, Japan, Kondo entered the Aeronautical Laboratory, became a lecturer, an associate professor and, in February 1946, (via appointments at Kyushu and Nagoya Imperial Universities) a full professor. He took charge of the 5th Chair of Applied Mathematics of the First Faculty of Engineering, Tokyo Imperial University. After a re-organisation of the university after the war, he became professor and Chair of the Department of Mathematical Engineering and Instrumentation Physics until he retired from the university in March 1971 [Sugihara, 2001]. He supervised the doctorates of a small number of prominent Japanese engineers and mathematicians.

After Kondo retired, retaining emeritus status, he continued to work more or less continuously at the Centre for Prevenient (i.e. Anticipatory) Natural Philosophy (CPNP), based at his home in Yotsukaido. He was still active academically within six weeks of his death. Two papers were published post-humously confirming the existence of a small infrastructure of up to a dozen individuals (former colleagues and students) that supported his work. The author briefly met, in London, Professor Shun-ichi Amari, one of his doctoral students and supporters, also a professor of mathematics at the University of Tokyo. Kondo suffered several year long periods of ill-health, which merely interrupted his academic output.

Kondo travelled little due to ill health. He made two trips to Europe in the 1960's and was the Chair and Chief guest of the Einstein Centenary Symposium in Nagpur, India in 1980. He spent a month or so in Sweden via Moscow in 1981. He found these trips physically exhausting. The proceedings of the Nagpur symposium [Kondo, K. & Karade, T.M., 1980] provided us with the first pictures of him, his wife and some colleagues. See Figures 1 and 2. We reproduce his Inaugral Address to the Nagpur Symposium in Appendix C.

Kondo, colleagues and other people interested in his work met regularly at the CPNP on the 2nd of January (Kondo's birthday) of every year with occasional later meetings in the summer. There were 47 such meetings which were documented annually in the Post-RAAG reports (see section 6). 15 people attended the January 1996 meeting. Foreign scientific visitors and their wives occasionally stayed over at the CPNP with Prof. and Mrs Kondo.

Kondo was a deeply spiritual person, which we can readily deduce from the occasional Buddhist or Shintoist footnote in his Post-RAAG work and by the curious symbol of Figure 3 which is embossed onto the front and back covers of the Research Association for Applied Geometry (RAAG) memoirs [Kondo, K., 1955, 1958, 1962, 1968]. This symbol also appears on the back cover of his last but one Post-RAAG monograph. Note that, for brevity, the post-RAAG monographs are referenced herein as [K001-K358]. The introduction to the first volume of the RAAG memoirs contains a translation of the famous Platonian slogan in Greek around the outer circumference "*God ever Geometrises*".

The significance of the four inner characters has not yet been deciphered. The general form of the diagram closely matches an electrical motor configuration diagram which appears in Volume I (p209) of the RAAG memoirs.

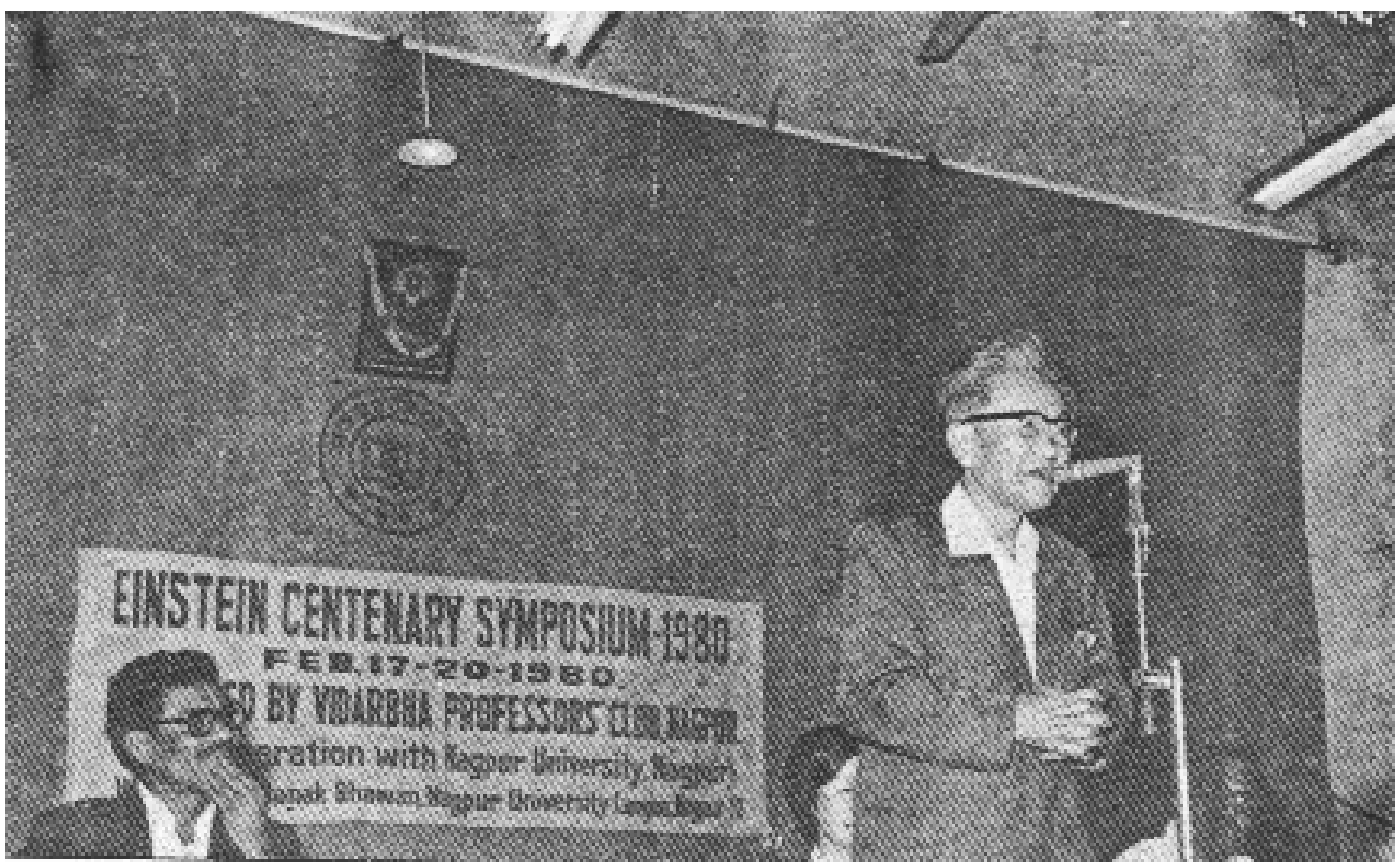

**Figure 1. Kondo making the opening address at the Einstein Centenary Symposium Nagpur, India, 1980**

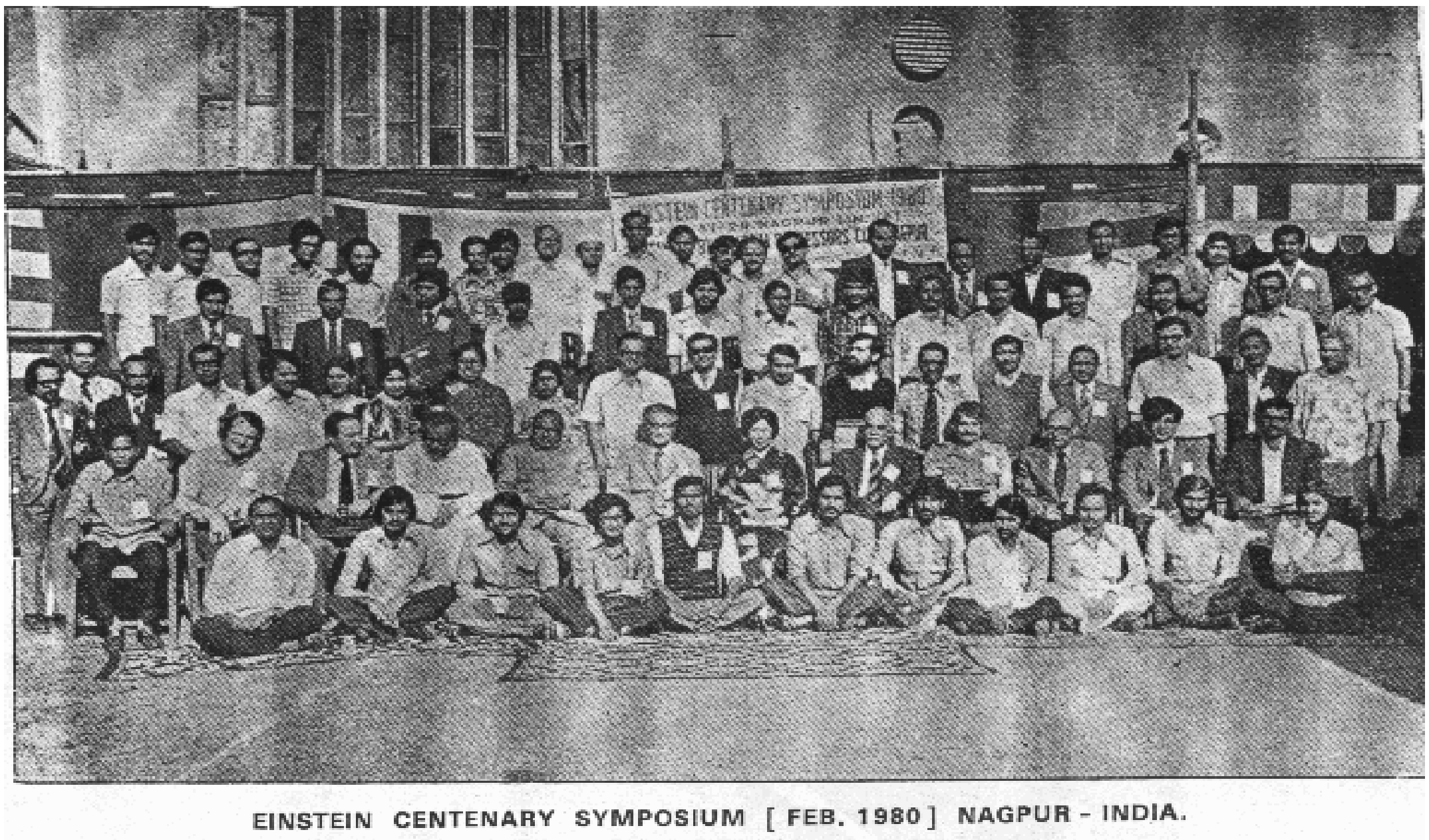

**Figure 2. Professor & Mrs Kondo with Symposium Delegates**

In [K030], Kondo writes of the transition from RAAG to post-RAAG:

> *"The RAAG appeared once to be a sanctuary where only true worshippers of God dared to enter. With an external prosperity it became gradually impregnated with paganism. Those who had once been angels having fallen, the Paradise was about to be lost, the spirit of the founding fading. A fire was lit in May 1973, for the angel's wingbeats to revive from the vestige as a phoenix plumed for a reparation campaign"*

It must be made very clear at the outset that Kondo's spirituality did not directly interfere with his work, though it clearly guided him. He had a dry sense of humour too, as we might also judge from the occasional footnote. Kondo was fluent in English and Japanese and had a varied and often impressive written command of other languages and character sets.

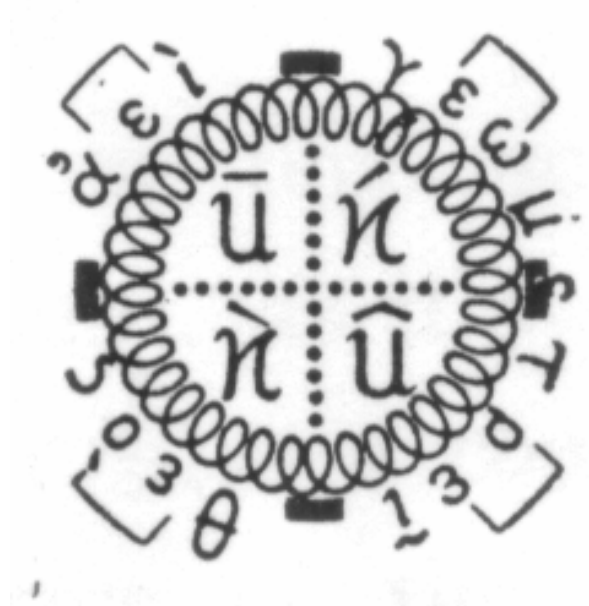

**Figure 3. Symbol embossed on the front and back covers of each volume of the RAAG memoirs**

Kondo was an honorary foreign member of the US National Academy of Sciences and an Honorary member of the Tensor Society of Great Britain. Kondo was elected on the 20th December 1973 a non-resident full member of the Neapolitan learned society Accademia Pontaniana. At the election, emphasis was placed on Kondo's contributions to phonetics, musicology and Prosody [K013]. He was awarded the Mainichi Academic Encouragement award of 1952.

This biographical information regarding Professor Kondo was kindly provided by: Tomoaki Kawaguchi, a professor of mathematics at Mejiro University, Japan [Kawaguchi T., 2004]; Professor Kojiro Kobayashi, SOKA University, Tokyo and Professor Laxmi Chandra Jain of Jabalpur, India. Copies of his 40 year correspondence with Professor Jain has also been graciously provided. See the Supplementary Materials section for an index and information regarding access to these sources.

## 2 Early Research Work

### 2.1 Aeronautical Engineering

Kondo graduated from the Department of Aeronautical Engineering, Tokyo Imperial University in the early 1930's. He then published some related early work in the mid 1930's [Kondo, K., 1935] which was well received. This work was interrupted by the Second World War and he was reluctantly forced to work in other academic domains..

### 2.1 Research Association of Applied Geometry (RAAG)

Kondo established in 1954 and then organised until 1973, the Research Association for Applied Geometry (RAAG). The purpose of RAAG

> *"…is to further the essential advancement of Engineering and physical sciences by the Unifying Researches of the Basic Problems through the medium of Geometry and other Mathematical Thought and Training*" [RAAG, 1971]

RAAG started as a 3 year government funded research project. There were initially 16 collaborators. They published an irregular series of research notes in the period 1951-1973 outlining their initial thoughts. The RAAG research notes [Oyo Kikaqaku, 1951-1973], amounting to 200 in all, are rather untidy and were used as an informal communication medium. The print run amounted to 280 copies suggesting that the original group grew, and that there was academic interest in the work outside of the group itself. The RAAG research notes ceased upon Kondo's retirement from the University of Tokyo in 1971.

Kondo was an early contributor to Matrix and Tensor Quarterly [Kondo, K., Amari, S., 1963] where he met or came into contact with its British founding editor Austen Stigant, a retired consulting engineer. Stigant provided Kondo with valuable assistance in the editing of the RAAG memoirs.

The 1971 RAAG list of members shows that there were only ever 276 members of RAAG (of whom Kondo is listed as member number 1). In 1971 there were two honorary members, Professor Akitsugu Kawaguchi and Mr Austen Stigant plus two Patron members Nishi Trading Company Ltd. and Kanto Auto Works Ltd. 204 individuals were listed as current members, suggesting an organisation in good membership health.

The RAAG research notes were substantially improved, collected together, peer reviewed and published commercially as the RAAG Memoirs in four handsome volumes during the period 1955-1968. They are occasionally advertised on Internet Bookshops – the author being the owner of A.J. Good's (RAAG member 144) copy of Volume IV. There are 116 papers in 2500 pages covering initially 5 research areas, later expanding to 10. All four volumes were edited by Kondo and 53 (45%) of the papers were also authored by Kondo. Publication was hampered by financial and other difficulties, including a year long period of illness suffered by Kondo. The RAAG research notes continue for a few years after the publication of Volume IV. The latter papers are indexed in [K008].

The initial RAAG research Divisions were:

Linear Geometry and Topology of Networks
Differential Geometry of Engineering Dynamical Systems
Geometry of Deformations and Stresses
Non-Holonomic Geometry of Plasticity and Yielding
Miscellaneous Subjects

The later Research Divisions added:

Geometry of Observation and Structurology
Diakoptics (Tearing, Tensors, Topological Models)
Topological Information Theory of Engineering System Structures
Mathematical Foundations of Psychophysical Recognitions
Geometrical Theory of Natural Information Patterns
Miscellaneous Subjects

Contributors to the RAAG memoirs included Gabriel Kron, Otto F. Fischer, Shun-ichi Amari, Mikiaki Kawaguchi and Tomoaki Kawaguchi.

Kondo's most significant early RAAG work was in Elasticity, Yielding and Dislocation. He discovered in 1952 the relation between dislocation density and Cartan's torsion using differential geometry and noted that the propogation of dislocation in crystals is relativistic [Kondo, K., 1955a]. This work was progressed much further [Bilby, B. A., Bullough, R., and Smith, E. (1955)], [Kroener, E. 1958, 1980].

Otto F. Fischer, a Swedish consulting engineer, published his unification ideas in the RAAG memoirs, supported by Kondo. Kondo was clearly inspired by them and took the ideas much further, referencing them as late as 1997 in [K316].

There are a number of papers by Kondo in the RAAG series which one might regard as key. The first being a paper entitled "Quantum Hydrodynamical Analogy by Isotropic Observation of Turbulence" [Kondo, 1958, p303]. The paper shows how Schrödingers wave equation appears in a theory of turbulence and includes an equivalent of Planck's constant.

> "*It is surprising that these scalar equations obtained from the purely hydrodynamical criteria by natural isotropifications, have the same form as the Schrödinger wave equation of a single particle in an electromagnetic field*".

It is a recurring feature of Kondo's work that he shows how the same mathematical representations occur in very widely differing fields.

Some Kondo papers are enigmatic: "Construction of Kawaguchi Space-Time by Statistical Observation of Monads and the Origin of the Quantum, Physical Fields and Particles" appearing in the RAAG memoirs volume III. This is Kondo's first unification paper, and also the first to be based upon Kawaguchi's mathematics. It follows Otto Fischer's quaternionic unification paper in the previous volume [Fischer O., 1958]. Within this paper, a footnote states:

> "*The argument of this subsection may not essentially be alien to another logical algebraical approach made by E.W. Bastin and C. W. Kilmister*",

There follows a reference to their first Concept of Order paper [Bastin E.W. and Kilmister C.W., 1954]. Bastin [Times, 2011] and Kilmister were founder members of the UK based Alternative Natural Philosophy Association (ANPA) of which this author has been a member for 15 years. In addition, there were close links between ANPA, RAAG and the Tensor Society of Great Britain whose President and founding editor was Mr. S. Austen Stigant.

There are also joint papers by Kondo and Michiaki Kawaguchi [1968] and Kondo and Tomoaki Kawaguchi [1968], perhaps signalling deep collaboration and relevance between the interests of Kondo and the work of Akitsugu Kawaguchi.

A deeply respectful obituary of Gabriel Kron appears in Volume IV of the RAAG memoirs:

> "*….the world of applied geometry loses the most outstanding pioneer*".

**3 Akitsugu Kawaguchi (1902-1980)**

It would be inappropriate in a biographical paper on Kondo's work to omit mention of the life and times of Akitsugu Kawaguchi. Kawaguchi was Professor of Mathematics at Hokkaido University, Japan. Founder of the Tensor Society, he established in the 1930's the Japanese school of thought developing the Geometry of Higher Order Spaces. The work was published in the Proceedings of the Imperial Academy of Japan. Unlike Kondo, Kawaguchi was well travelled, and as with Kondo, held in very high international regard

> *"...Not only in the field of pure scientific research, but also his contribution to mathematical education in schools in Japan is so great that one cannot disregard it*" [Kawaguchi T., 2003].

Kawaguchi had 7 children. Four sons became professors of mathematics. Tomoaki Kawaguchi was a recent president of the Tensor Society. Kawaguchi's loss was deeply felt by Kondo as evidenced by an obituary which later appeared. Of Kawaguchi's work Kondo writes:

"*The discovery of higher order spaces, construction and development of a fundamental theory of their geometry impress us with such a power as we recall them on this occasion*" [K179].

Kondo first adopted the Kawaguchi Space as his lingua franca within the RAAG memoirs of 1968, continuing to use them almost without interruption until his death in 2001.

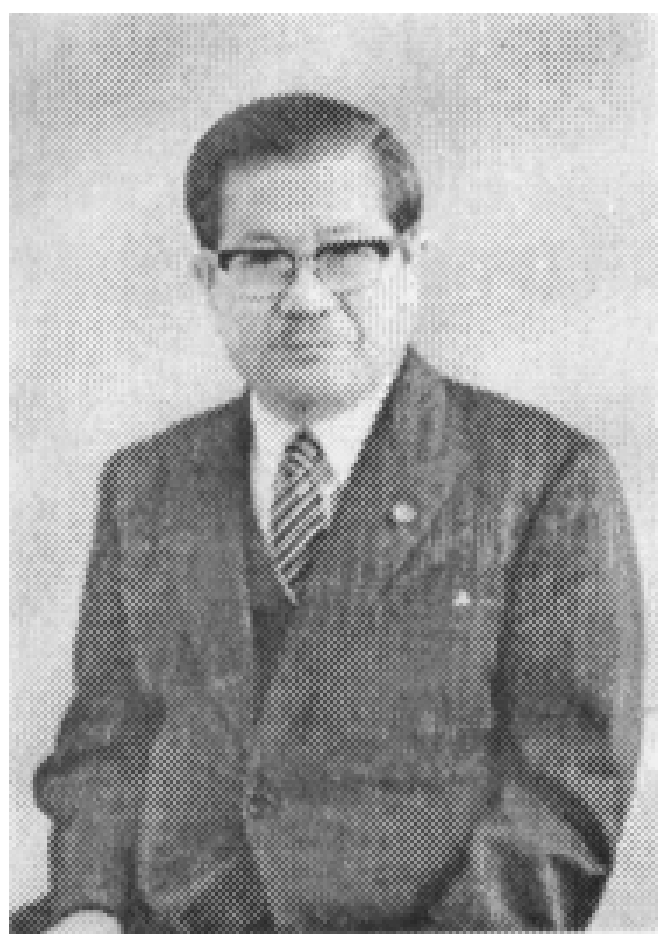

**Figure 4 Akitsugu Kawaguchi**

For those who are unfamiliar with Kawaguchi spaces [Kawaguchi A.,1931, 1932, 1933], there are two important introductory papers in RAAG Volumes III and IV by Michiaki Kawaguchi, comprehensively outlining his father's work [Kawaguchi M. 1962, 1968].

**4 The Kawaguchi Space**

The Kawaguchi space was first so named by J.L. Synge [1935]. Very briefly:

The Riemannian arc length of the curve *k(t)* is given by the integral of the square root of a quadratic form in *x*' with coefficients dependent in *x'*. This integrand is homogenous of the first order in *x'*. If we drop the quadratic property and retain the homogeneity, then we obtain the Finsler geometry. Kawaguchi geometry supposes that the integrand depends upon the higher derivatives *x''* up to the *k*-th derivative $x^k$. The notation that Kondo uses is:

$$\mathbf{K}^{(M)}{}_{L,\,N}$$

For:

| | |
|---|---|
| L | Parameters |
| N | Dimensions |
| M | Derivatives |

Much is made of the use of Kawaguchi Trees (Figure 5) and its three pronged form. The extended branches in Fig. 5 represent the higher orders, M, as numbered.

Within the realm of sub-atomic particles, Kondo makes very extensive use of Kawaguchi Spaces, and their mathematical properties are repeatedly demonstrated to be equivalent to the observed form and measured performance of classically understood physical entities.

Numerous Kondo papers cover the equivalence of the spaces $\mathbf{K}^{(M)}{}_{L,N}$ and elementary particles.

“*Baryons are identified for L=3, Mesons for L=2 and Leptons for L=1. The physical configurations carried by the three / two parameters of* $\boldsymbol{K}^{(M)}{}_{3,\,N}$ *and* $\boldsymbol{K}^{(M)}{}_{2,\,N}$ *are identified with the three / two quarks of the baryon / meson*”.

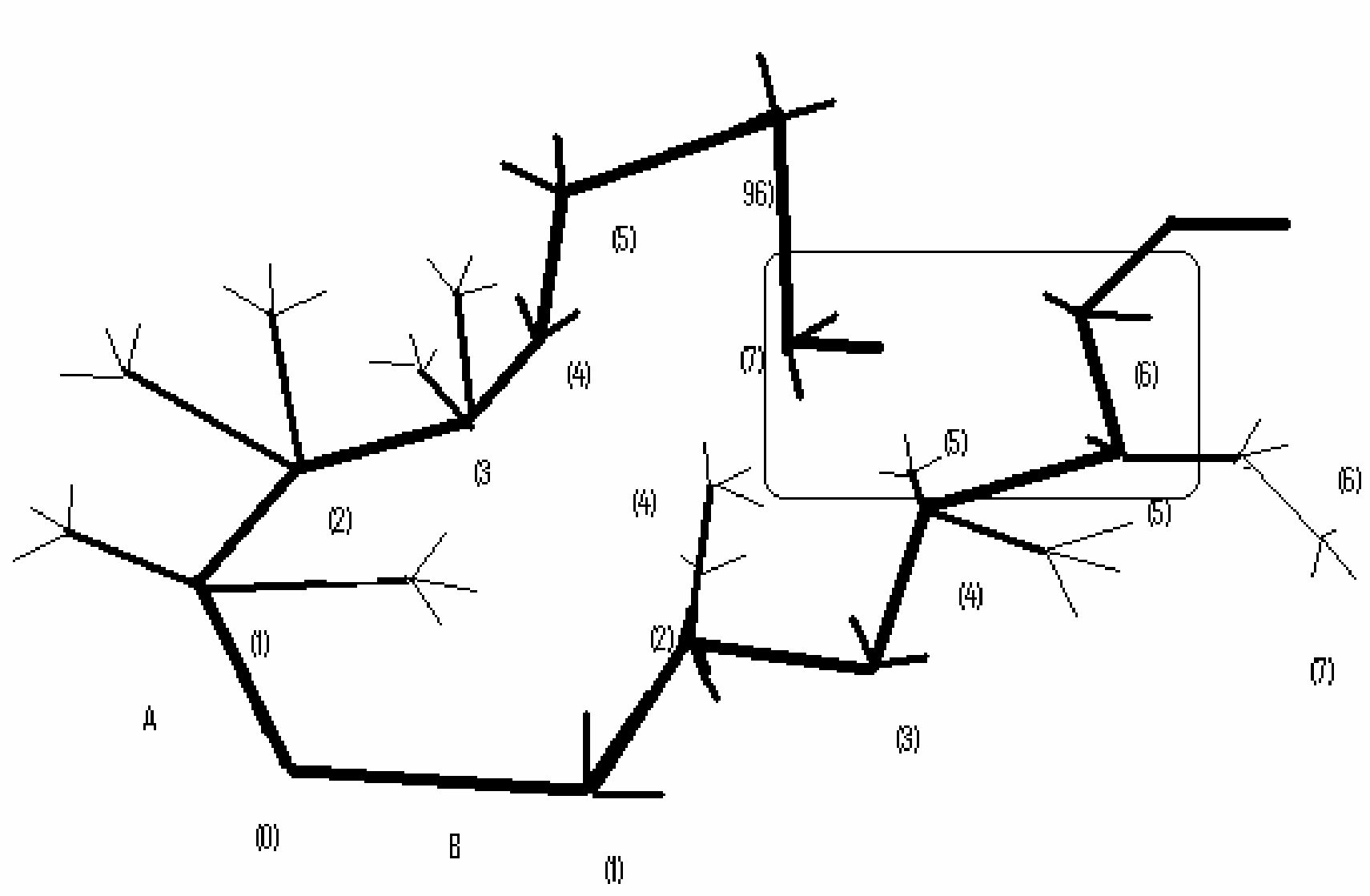

**Figure 5 A Kawaguchi Tree**

A further result, due to Kondo, the KHK theorem [K291, p21], appears to provide a mathematical basis for the three dimensional nature of our observed universe. The evidence for this is provided not only by the result itself and the context, but through the almost monotonous repetition of the theorems, in summary and detailed form throughout the 28 year period of the Post RAAG reports.

## 6 The Post RAAG reports

Following his “retirement” from the University of Tokyo, Kondo set about continuing the work of RAAG, through a series of generally monthly Post RAAG reports. Entitled “Research notes and Memoranda of Applied Geometry for Prevenient Natural Philosophy”, the series is a breathtaking and large piece of research, covering a very diverse range of topics in the nuclear, physical, chemical, biological and other natural sciences, including the study of communication, language, linguistics, phonetics, acoustics, musicology, prosody, consciousness and perception. It was these notes, lying on a shelf in the British Library that first drew the author’s attention, in 2001, to this vast body of work. The author had been studying the subject of cryptography [Nichols R.K., 1998], the journal Cryptologia being shelved “in close proximity” to Kondo’s Post RAAG reports. The 10 year shelving policy at the British Library means that these two bodies of work are no longer openly shelved together and now lie separately in deep storage. See Appendix B for the outline of a three day lecture given at Hokkaido University in 1976 which maps out Kondo’s strategic perspective.

In total there were 358 Post RAAG monographs published from 1973 until Kondo’s death in Dec 2001. There were very few temporal gaps. Each monograph comprised 8-100 pages, resulting in approximately 7,500 (B5 – Japanese industrial size) pages total. The monographs were almost entirely in English, with a one or two page Japanese summary. The monographs were professionally typeset to a high standard, and the mathematics was mixed with substantial textual explanation, diagrams and tables, making them accessible to the curious non-specialist. The index is too large to be reproduced here, but is available in the Supplementary Materials. It is hoped that the Post RAAG reports will continue to be scanned and further placed in the Supplementary Materials. [K333] provides us with a colour picture in Figure 6 of Kondo on his 88th birthday.

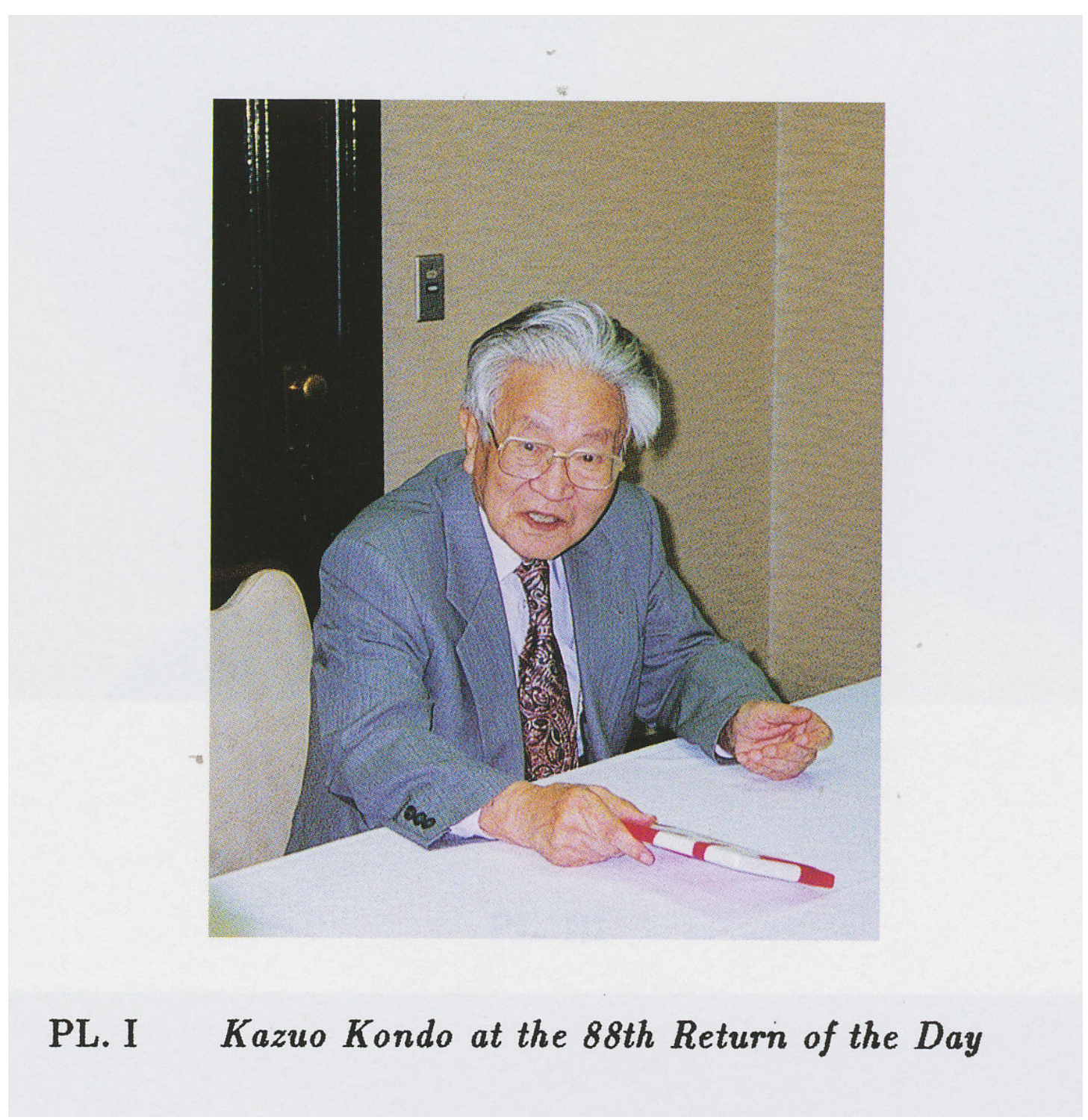

PL. I *Kazuo Kondo at the 88th Return of the Day*

**Figure 6 Kondo on his 88th birthday**

Kondo authored about 90% of the Post RAAG reports, however there were some very interesting contributions from Karade (on spherically symmetric space times), Thomas, Butey, Jain, Amari, Uehara, Gluckman and others. Each year an annual index was published which contained a brief report of the Jan 2nd and other meetings, additions to the CPNP library and the occasional obituary. There was an annual fee of 4,000 Yen. The Post RAAG monographs were donated to a handful of national libraries (British, Smithsonian etc), and are exceptionally rare. The author is grateful to Professor Tomoaki Kawaguchi for obtaining an almost complete set of Post RAAG reports for his personal use. The post-RAAG print run was 200 and there were up to "just under 100" subscribers.

The supplementary materials contain the entire correspondence from Prof. Laxmi Chandra Jain, of Jabalpur, India to Prof. Kondo over the period 1963 – 2001 regarding the publication of Jain & colleagues 20 papers for Post-RAAG and Kondo's visit to India in 1980. The correspondence includes a 2002 letter from Prof. Jain to Prof. Amari regarding Kondo's passing and the accidental copying of a small part of a handwritten letter in English from Mrs Reiko Kondo to Prof. Jain containing the poignant sentence *"I miss him very much"*.

The supplementary materials also contain a long administrative record of post RAAG activities. This material was kindly provided by Professor Kojiro Kobayashi. The correspondence includes the invitations to Kondo's 77th and 88th birthday celebrations.

During the early years of the Post RAAG reports, subjects were apparently studied at whim. The first report, K1, studied interference in Black and White television pictures and somewhat characteristically showed that the mathematics of such interference patterns was in common with the mathematics of elasticity. [K035] was on Quaternionian Phonology. [K068] was on the Life Cycle of Growth, Metamorphosis, Ecdysis, Mutation, Zygosis, Dissemination, Procreation and Death. [K137] was, naturally, devoted exclusively to the fine structure constant where an *a priori* derivation of $5\alpha^2 + 137\alpha - 1 = 0$ is given, ie $1 / \alpha = 137.036(4793)$.

A recent reading of [K-46] of March 1976 suggests that Kondo's approach was in fact highly structured. Kondo gave a three day lecture series on February 19 – 21 1976 at Hokkaido University outlining his research agenda. We reproduce the agenda in Appendix B. There is a very strong link between the subjects covered in the 1976 lecture series and the contents of [K001-K358]. It is interesting to note that two research topics were "*Geological and celestial-planetal organisms*" and "*Organic Cosmos*" [Lovelock, 1979].

## 7 Major Post-RAAG Research Series

The later work was organised into a number of series on major themes, however study of the various themes would temporally overlap. The last series, a complete re-organisation, was started in 2001 but was unfinished. Here follows a by no means complete listing of the areas covered by these later major series.

### 7.1 Epistemological Foundations of Quasi-Microscopic Phenomena from the Standpoint of Finsler's and Kawaguchi's Higher Order geometry (Jun 1991)

The papers in this series are: Introduction, Quantum and Space Construction; Electromagnetism; Gravitation; Meaning of Defectiveness; Fourth Order Physics; Yielding; BRIEF with a prospect for furtherance

### 7.2 Higher Order Geometrical Theory of Elementary Particles (May 1993)

The papers in this series are: Introduction and Theory of Low-Lying Leptons; Theory of Hadronic Multiplets; Higher Flavours; The Theory of Massive Gauge Bosons; Summary

### 7.3 Strange Stories of the Mathematical Constitution of the Bio world (Oct 93)

The papers in this series are: Introduction and Perspective; Standpoint of biocreatures; Morphology of Vascular Plants, Classification of Vascular Plants; Floral Metamorphosis; Metric Aspects of Leaves and Roots; Phylotaxis and Inflorescence; Zoological Phase, Locomotive devices of Animals; Evolution of Tetrapodal Vertebrates; Animals with radial symmetry construction; Uncertain Aspects of Bio-Elementary Lower Bounds, Bio Feasibility at Cytological level; Physio-Psychological Information.

### 7.4 A Unified System of Epistemological Penetration to Nature (Oct 1995)

The papers in this series are: Introduction and Fundamentals [K291]; Meaning of Quantum Mechanics; Construction of Elementary Particles; Atoms as Hyper-Multiparametric Space; Inter-Molecular Organic Chemistry; Topological Construction of Microscopic Observations; moving and Standing Particles; Polyhedral Aspects of Atoms and Molecules; Order and Dimensional Preference by the Generalised KH and KHK theorems; Yielding as a Catastrophic Deviation from the GKHK Limit;

It is by reading this major series of monographs that one first comes to realize that one is reading something very special indeed. The paper entitled "Meaning of Quantum Mechanics" (K292) gives a clear reconciliation of the geometry of higher order spaces to quantum mechanics as conventionally understood. The paper entitled "Construction of Elementary Particles" (K293) gives an *a priori* Higher Order space derivation of the Weinberg Angle:

$$\text{Cos}\ \theta_w = 17 / 19 = 0.8994737$$

And notes that "*The Fermi constant is but a transform of the Higgs mass.*.". The next section outlines the first paper in the above series.

### 7.4.1 Introduction and Fundamentals (K291)

In repeated readings of Kondo's work, several papers appear particularly prominent, of which K291 would be of most relevance to a wider audience. We reproduce verbatim what seem to be some pertinent paragraphs in the introduction to this paper and leave the mathematics in the body of the paper to the specialist, hoping that this introductory material will inspire the reader's curiosity. Kondo writes:

> *"The author's investigation has been motivated even more than thirty years since with a discovery of a resemblance of the relativistic wave equation for the elementary particle, especially for the electron, to the fundamental condition equation for the geometry of higher order spaces. Pursuing the clue so-suggested, he has arrived at the recognition of an epistemological thesis:* ***Higher order spaces are themselves particle and particles themselves are spaces****".*

He continues:

> "*Now, the space used to be treated as a continuum, described in terms of coordinates which are real numbers. However, we should not forget that real numbers are products of a compromise to adopt among numbers limits of sequences of rational numbers. On the other hand, the most basic recognition starts with fundamental information elements which should not have individuality and can only be counted. The natural numbers obtained by counting are amplified to the ring of integers and then, by comparison processes, to include fractions, to the field of rationals, toward the contiguous field of real numbers. Thus we have connected the basically finite discrete microscopic edifice to the apparently continuous macroscopic world picture*".

In a couple of paragraphs, Kondo outlines his world view, suggesting that the basic unit of the universe is an indistinguishable information unit, that particles are merely mathematical spaces and vice versa, giving an immediate dualistic connotation. In the following paragraph, Kondo outlines the boundary mechanism between the micro and macro worlds:

> "*As is well known, one can start with dividing the macroscopically observed object occupying ordinary space into cells or 3-dimensional simplexes by triangulation. Each such cell is bounded by a closed surface. The boundary can be further triangulated into area elements or 2-dimensional simplexes. One can further proceed to a decomposition into 1-dimensional simplexes or line segments, each bounded by a pair of points or two 0-dimensional simplexes. The decomposition need[s] to terminate at this step macroscopically. However the mathematician considers further advances to negative dimensional simplexes of -1,-2,….,-∞ dimensions or what one calls dual cells or cocells or dual simplexes. Pertaining by definition to regimes macroscopically unobserved, these cannot but be of the microscopic world. Each cell and/or co cell needs to carry some phase of macroscopic/microscopic physics*"

Thus the richness and variety seen within sub atomic particle physics is a manifestation of the diversity of negative dimensional Kawaguchi spaces of various orders. In other papers [K293] we see the tetrahedron as being the basic information carrying unit of the nucleon, with information carried on the faces, vertices and edges.

Kondo continues:

> "*Cocells formally multiplied by microphysics representing coefficients form a linear algebra with the respective dimensional standpoints. The mathematician calls the elements of the algebra algebraic complexes or chains or dual complexes or cochains. As a boundary / coboundary is associated with each cell/cocell there is also entailed an* *algebra of boundary / coboundary complexes. The operation to derive boundaries from chains is ∂ and the operation to derive coboundaries from chains is δ. A most fascinating and important aspect of these operations is that applying either ∂ or δ twice always yields an empty algebra: ∂∂ = 0, δδ= 0. Since the repeated application sequence of boundary / coboundary operation[s] terminate [at] once at this, the mathematical formulae in theoretical physics are described mostly in terms of second order differential equations. This is especially the case for the classical description of macroscopic physics, whether of particles or of continua. If more microscopic phases are penetrated and observed to be working in cooperation with the macroscopic phases, more higher order differential equations have to be handled. Therefore, differential equations of fourth order have to be introduced for the description of this sort of quasi-microscopic phenomena*"

"*Put in terms of expoint co-ordinates, involving not only simple coordinates of N dimensions, but their derivatives up to a finite order, say M, with respect to a finite number, say L of parameters, the geometry of higher order spaces* $\boldsymbol{K}^{(M)}_{L,\,N}$ *are to all intents and purposes in correspondence with the boundary / coboundary penetrations of microscopic spheres. Thus we may assume that we are prepared fundamentally to decipher the mysteries of microscopic physics*."

Thus we have reached page 3 of the introduction to K291 which in itself is an introduction to Kondo's epistemological world view. We continue with a few paragraphs outlining the body of the paper.

"*Chapter 1 starts with a description of this recognition and is concluded with a reference to classical mechanics as an illustration of its specific case restricted to the sphere of* $\boldsymbol{K}^{(1)}_{N}$ *of order 1. How the classical viewpoint is so restricted is illustrated with reference to the Hamiltonian principle, restriction to causality, contact transformation, so entailed*".

"*Quantum Mechanics dealt with in Chapter 2 are also not entirely free of the same restriction. Emphasis has, therefore, to be placed on the concept of Hamiltonian which has originated from the classical mechanical construction although some higher order considerations work for letting the wave-mechanical formalism to be preferred. Normally an osculation of higher orders to the order 1 is introduced so that quantum-mechanical routines used to be put without reference to the latent higher order expoint background. If the correspondence principle has not so far be[en] able to provide more than an ambiguous impression, the crime is that the exposition used to be without reference to this hidden mechanism*".

"*We reconstruct wave geometry, proving that the wave equations are but reinterpretations with some osculation approximations of the Zermelo conditions. However the Schrödinger equations comes only for*

$$\Delta_1\Delta_1 F = F$$

*By osculation this is brought to a partial differential equation of second order*".

"*But* $\Delta_1 F = F$ *being more basic, partial differential equations of first order must come first with less reliance on osculation. We shall argue for the relativistic wave equation for the electron from this point, to be reconciled to the spinorial language.*

*Our reinterpretation of a further characteristic point of quantum mechanics has also been prepared by the geometry. Kawaguchi has drawn attention to what he calls Modified Zermelo Conditions, relevant for the system of two higher order line elements stationed mutually very closely* [note the Box in Figure 5]. *The condition equations have*

*the structure that is not concerned with how many of such systems come together while the ordinary Zermelo conditions depend upon the number of line elements the geometer calls the weight. Two classes of particles are, therefore, spontaneously distinguished from each other, as follows.*

*An arbitrary number of particles under the modified Zermelo conditions can be carried on one single parameter, i.e., located within one and the same configuration. This is just the standpoint of bosons, leading to the Bose-Einstein statistics.*

*The same privilege is not accorded to particles under the simple Zermelo conditions. The construction of the conditions depending upon the weight, not more than a single particle of this class cannot be located in one and the same configuration. This is no doubt the standpoint of fermions, leading to the Fermi-Dirac statistics.*

*Further discussions in Chapter 2 are concerned with spin and electric charge. This latter is something first arising on reaching the order 2, which means for the electron / positron. It will be shown that spin ½ is carried by a single line element whichever the order. This is the primary form of the fermion. A boson of primary form has naturally spin* $0^-$ *or* $1^-$ *from our theoretical standpoint. It will be argued that fermions according to the quantum statistics can be of half integral spin and bosons can be of integral spin".*

*"Among important clues that have been obtained, emphasis can be placed on the point that the masses of particles are quantized into sequences of quanta of different sizes consecutively multiplied by a constant to form a scale:*

$$\ldots\upsilon \;\rightarrow\; \theta \;\rightarrow\; \Theta = \theta / \alpha \;\rightarrow\; ¶ = \Theta / \alpha \;\rightarrow\; \ldots$$

*Where* $\alpha = 1 / 137.036\ldots$

*is known as the fine structure constant. So far the existence of the natural constant having this value has been empirically well confirmed. However this is something which has to be proved theoretically……the author has had himself to attempt his own approach, with the help of higher order geometry, to clarify by what epistemological construction more minor quanta are eventually composed into the 1/ α value."*

K291 extends to 30 pages, followed by the nine continuation papers outlined in section 7.4 above. Note the introduction of a series of mass quanta.

### 7.5 Latent a priori Frame for Mammalian Human Morphology, Anatomy and Physiology (Apr 97 )

The papers of this series are: Introduction and Fundamentals; Neural-Cerebral Design; Internal Organs; Osteo Myological Quantization.

Again, note the use of a scale of mass-quanta as in section 7.4.1. In the paper on Osteo Myological Quantization, bones are viewed as a higher order quantization. Note also the use of the word design in the title of the second paper of this series. In order to illustrate the breadth of Kondo's work, we again quote liberally from the paper on Osteo Myological quantization, which deals with the geometrical origins of the skeletal features of living beings.

#### 7.5.1 Osteo Myological quantization (K258)

"*The site of the parameters in a higher order space can also be quantized into segments, the limits of which can be no more decomposed. Such a limit may be nearly a rigid piece. In the animal body such quanta cannot but be bone pieces forming parts of the skeleton,* *whether lying internally as [endo]-skeleton or as almost rigid shell covering the body as external skeleton*."

"*Note the partition of the body into three main segments: Head (cephalique), pectral (breast), caudal (tail), materializing the KH order limit M >= 3 or the KHK dimensional limit N >= 3. Notice also the quantization into more macroscopic segments such as of the abdominal part into several smaller segments beyond the KHK lower bound N = 3. Lateral symmetry with a symmetry axis is remarkable. This is of course an indispensable consequence of the modified Zermelo conditions, which entails also locomotive appendages differentiating into legs for walking and wings for flying in the case of insects*"

Two paragraphs addressing the simple issues of what bones are, mammalian bi-lateral symmetry, the numbers of major body parts and their segmentation, the notion of the mathematical origins of wings, legs and arms. The author recalls reading in other papers the dimensionality of eggs being zero, hence their need of warmth for progression to locomotion and the dimensionality of snakes being one, hence their mode of locomotion. A feature of the biological papers is their attention to detail, their use of line art to depict the various forms of living being – from birds to starfish to dinosaurs, the use of the full latin terminology and at all times the relationship of the various form of living being to the underlying higher order geometry and the mathematical notion of principle ideals. In Figure 7, the human skeleton is treated as a hierarchical Kawaguchi tree with its characteristic three pronged form.

*"The lower part of the skeleton can be divided into three prongs, each starting from the centre as a single parametric Kawaguchi tree*"

*"...the skeletal, muscular, gastrointestinal, circulation systems etc combine into a holo-parametric whole that can be more generally quantized, each quantum involving some osteological, neural, circulatory functions etc*

*"…thus globally the human body from head through trunk to limbs are quantized into a finite number of quanta*"

**7.6 Perception of Space, Matter and Bio-Psychological Information (Nov 97)**

The papers in this major series are: Galois Field, Multi Dimensional Space and Differential Geometry (K316); A Reorganisation of Mathematical Musicology; Confrontation of Linguistics and Musicological Sound Quantization; Undecidability and compromise Inevitable for the construction of musical scales; Golden Section Frame Restricting Visual Perception and Growth of Bio Creature.

This series addresses perception, consciousness. K316 is again outstanding, as within a few pages Kondo illuminates his main thesis regarding number and reality, outlines the nature of consciousness as "…*the extraction of Galois fields*", then progresses directly into mathematical musicology by way of illustration. Gödel's incompleteness theorem features prominently in the introduction to K316.

**7.7 A More Fundamental Analysis of Epistemological Penetration to Nature and Restrictions Inherently Latent in it (Sep 98)**

The papers in this major series are: Introduction, Combinatorial Topological Penetration into the Microscopic World; Problems of Heavy Elementary Particles; Constitution for Biological Classification; Definition and Stability of Elementary Particle, Neutrino Oscillation, Weak Interaction, Higgs Mechanism; Mass Contribution Spectrum of Quarks.

This series is his last rework of his main ideas in theoretical physics.

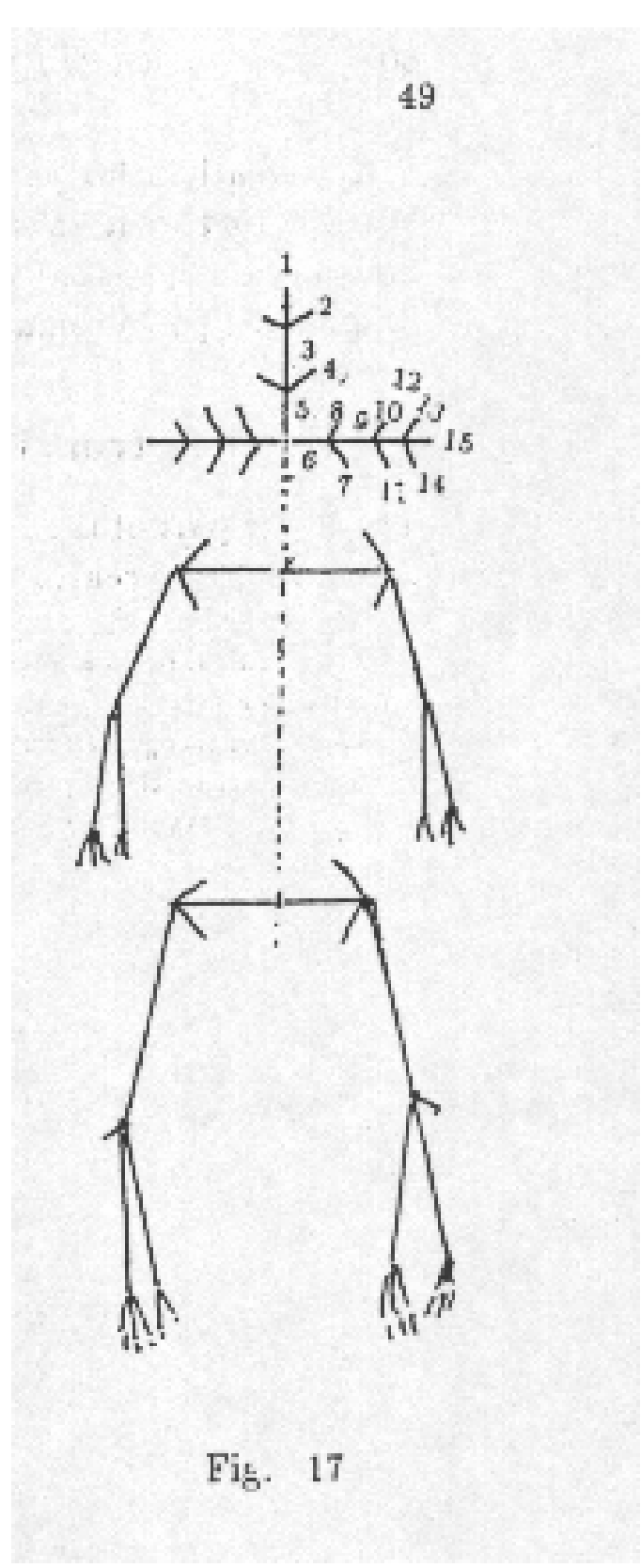

**Figure 7. "Decomposition of Holo-parametric skeleton into apparently single parametric units"**

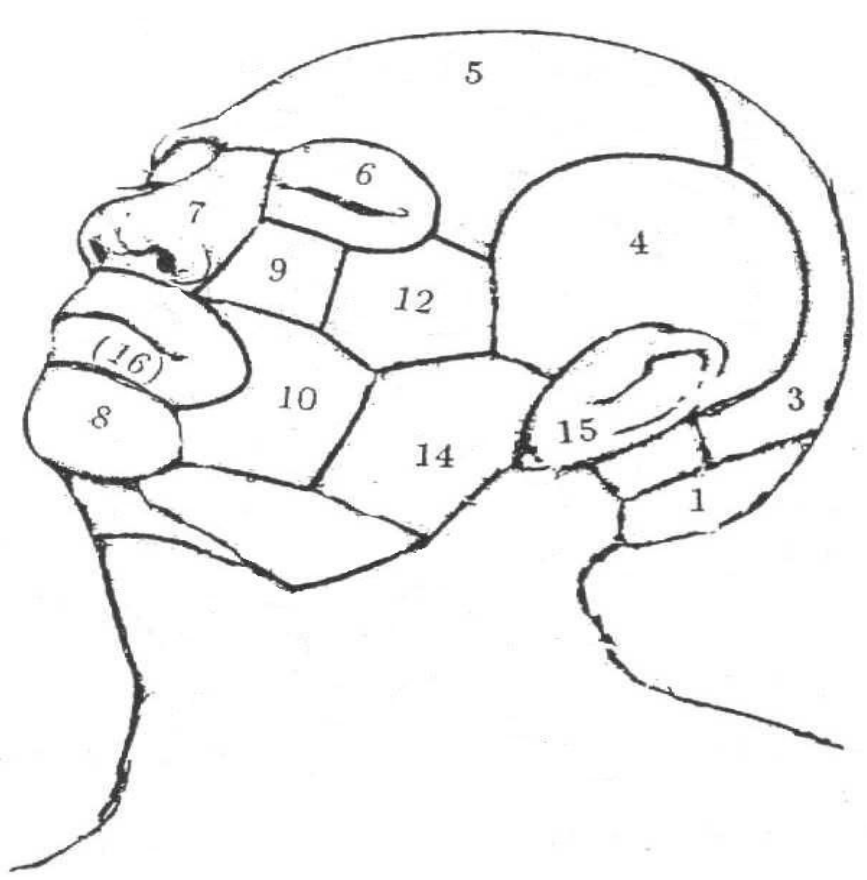

**Figure 8. "Quantum Sub Systems of the Neural-Cerebral System"**

"*The central single-parametric branch of the head carry regions 1-5 and the bilateral part is subdivided into paired regions 6-15*"

1. *R. Sphenoidale*
2. *Regio occipitalis*
3. *R. Paritalis*
4. *R. Temporalis*
5. *R Frontalis*
6. *R. Orbitalis*
7. *R Nasalis*
8. *R Mentalis*
9. *R Infraorbitalis*
10. *R Buccalis*
11. *Maxilla*
12. *R. Zygomatica*
13. *R. Mandibula*
14. *R Masseterica, oralis*
15. *R. Oralis*

**7.8 Re-organisation of the Mathematical Constitution of Natural Philosophy towards the New Millennium (Apr-May-Jun-Jul 2001)**

This (K357) is Kondo's last paper prior to his death in December 2001. It is a lengthy monograph written entirely in Japanese "*to facilitate the flow of ideas*", save for the title, and a final poignant note of farewell. It was published posthumously. The mathematical symbology within suggests that he is recasting by now familiar ideas. The back cover has the curious symbol of Figure 3 reproduced upon it. A second post-humous paper (K358) appeared, with a note that all post RAAG activity had ceased.

## 8 Accademia Pontaniana

Kondo published two long monographs of his work in the Quaderni dell'Accademia Pontaniana [Kondo K., 1989, 1997]. This is one of the few outlets of his work in the Western Hemisphere. At the time of writing, both monographs are still in print.

The first monograph entitled "Tristimulus Foundation and Wild Topology of the Construction of the Phonetical Manifold" is a summary of Kondo's phonetical work over the previous twenty years. The mechanism of language and communication is mapped onto a phonetical tetrahedron. The mechanism is "*…..put in terms of the Pauli matrices leading to a quaternionian representation*". There is a section on "*The Eightfold way of semivowels*".

The second monograph, entitled "Three Phases of Epistemological Penetration to Nature", is a fairly comprehensive, but necessarily terse 100 page coverage of the main themes of Kondo's work over the previous 50 years, some of which has been covered earlier in this paper. The concluding remark is as follows

> *"The twentieth century has indeed been a century of pragmatical investigation having abandoned the question of 'why' exploring mostly 'how' profitably models can be made use of for direct human benefits or desires. [I] hope the next century be otherwise…".*

## 9 Main Observations & Conclusions

It is impossible, even in a paper of this length to give a proper sense of Kondo's life's work. To do so fully and properly will take many years. In addition, it is almost impossible to closely read all of Kondo's 10,000 pages. However when reading his work, numerous interesting quotations come to light which serve to give an insight into his thinking. We list some of these quotations in Appendix A.

In conclusion there are some clear threads that run through Kondo's work. Kondo clearly believes in an *a priori* schema – by mathematical necessity, or even design.

In common with much ANPA work, there is a realisation that the fundamental building blocks of the world around us are of the most simple mathematical kind – integers, monoids, indistinguishables etc.

Throughout his entire work, Kondo enjoys finding out the mathematical commonplaces amongst the most diverse and abstruse concepts. Throughout the latter fifty years of his life, he has focussed largely, but not exclusively, on the applicability of the geometry of higher order spaces, specifically Kawaguchi spaces, to explain a surprising number of phenomena. Whether these are the right explanations, only time and experiment can tell. Latterly, he shows an equivalence between the mathematics of Hamilton's quaternions and Kawaguchi's geometry [K316, section 1.5], which may well be of interest to Rowlands [Rowlands P.,1999].

Kondo's world is very clearly digital and bounded – as for example consciousness is simply the extraction of Galois fields. Kondo's world is a set of hierarchical monoidal data structures, there being a mapping through the mechanism of the Kawaguchi space and by virtue of the KHK theorem into our observed largely continuous observed three dimensional reality. The worlds of fractals and cellular automata are also accommodated within Kawaguchi's geometry.

There is a hierarchy of mass quanta – from the particle realisations at the smallest scale, through to bones at a human scale, the ratio of the smaller and next larger mass quanta being the fine structure constant. Amongst the operation of these quanta the mathematics of the Kawaguchi spaces apply at all scales.

The dualism we experience is reconciled through the equivalence of spaces and particles - through the reconciliation of a world of discrete indistinguishables with infinity and through the accommodation of Gödel's incompleteness.

Kondo's work is reconciled to the accepted theories of mathematical physics, but it seems to be *a better* theory, as the same mathematics seems to be applying across widely varying fields. The use of higher order differential geometry, as opposed to the uniform use of second order differential geometry seems to this relatively unspecialised observer, to be a defining differentiator.

To suggest that Kondo has achieved a unified theory of our universe *might* be going too far. It would certainly be imprudent for an inclusive scientific community to simply *ignore* this vast piece of work.

## 10 Supplementary Materials

A significant amount of biographical material, an index of Kondo's post RAAG work and scans of (at least some of) K001-K358 are available in the Figshare Data Depository: Croll, Grenville (2020): Kazuo Kondo. figshare. Dataset. https://doi.org/10.6084/m9.figshare.12413090.v1
The post RAAG reports are being gradually uploaded into Figshare in approximately reverse chronological order commencing with K291 – K301and K333.

## 11 Acknowledgements

The author once again thanks Professors Tomoaki Kawaguchi, Kojiro Kobayashi and Laxmi Chandra Jain for providing substantial reference, biographic, photographic and other materials.

**Appendix A - Kondo Quotations**

*"All particles are spaces, All spaces are particle"*

*"Every living creature, human, plant or animal, is also a kind of higher order space"*

*"The most basic perception starts with fundamental information elements which do not have individuality and can only be counted"*

*"The counted numbers cannot but be finite"*

*"The epistemological recognition, or consciousness, of human beings, must primarily start with the extraction of Galois fields"*

*"As the information space is finite and discrete, no accuracy beyond counting unity, revealed as the Planck Constant, can be achieved. Heisenberg's uncertainty is, therefore, natural and indispensable"*

*"Three dimensions is the lower mathematical bound for the information space to support sufficiently complex perceptions"*

*"Hence macroscopically observed phenomena are three-dimensional whether they are of physics, chemistry or biology. One can attribute to this mathematical lemma the reason why our observed space has a three dimensional appearance."*

*"The ultimate microscopic mass assembly must be of three dimensional form and of the simplest kind. It cannot but be a three-dimensional simplex having tetrahedral form"*

*"…we shall concentrate upon simpler problems concerning elementary particles demonstrating the equivalence of mathematical spaces and elementary particles and their classification into leptons, mesons and baryons for 1, 2 and 3 parameter spaces respectively"*

*"137 smaller mass quanta are assembled into a unit of larger mass quantum, or the smaller and large mass quanta occur in the ratio 1/137 which we name a [inverse] fine structure constant"*

*"It has been pointed out by mathematician Kurt Gödel that, whichever system of self evident axioms one may start with, one has to introduce something not apparently self-evident in order to account for certain epistemological indispensability…As mentioned by Gödel, the most conspicuous of deviations from self-evidence is the concept of infinity."*

**Appendix B – Hokkaido Lecture Series**

Contents of Professor Kazuo Kondo's Lecture Series at Hokkaido University February 19-21, 1976, on the "Theoretical Fundamentals of the Structurology of Human Recognitions of Nature": A Prolegomena to PREVENIENT NATURAL PHILOSOPHY

Part I. Fundamental Principles and the Physical World

Introduction

Implication of natural philosophy; Range of the positivism; Various mythological alters; Subject matters of natural epistemology

I. Finiteness and Unitarity

Finite recognitions and real co-ordinates; Unitary representations; Four-dimension pseudo-Euclidean space; The restriction to quantum mechanics; Unitarity of the world of elementary particles

II Information-Theoretical Standpoint, Statistical Approach and Entropy

Entropy as information amount; Inorganic physical world; Wave mechanics; Organic world; Psychological informations; Triality of the epistemological world

III Higher Order Geometrical Picture of the Construction of Space-Time

Simplest Unitarity; Conformal invariance; Higher order space picture; Three-dimensional space constitution; Space-time constitution; Principle of special and general relativity

IV Standpoint of Elementary Particles

Absolutely stable particles; Atmospheric fields of particles; Polarity of Leptons; Dirac's wave equation; $\pi$ and $\mu$; Hadron's mass spectrum

V World of Atoms and Elements

Inert Elements; Mod-3-singularity; Differentiation of metals and non-metals; Natural radioactivity

Part II. Biological and Physiological World

VI A Revival of the Concept of Urorganism

Definition of living things; Uroganism and base space; Ur-epiderm; Flow of matters; Flow of informations

VII Cellular Processes

Topological aspects of cells; Natural metamorphosis and fission; fusion and mating; Pair of projective planes

VIII Materials of the Biological World

World of Bosons; Upper limit of bio-elements; Inhibitions and bio-activity; Carbon-based organism; Protein and nucleic acid; Geological and celestial-planetal organisms ; Organic cosmos

IX Genetic Informations and Prosody

Molecular processes; Informations carried by side chains; Constructions of nucleic acids; space and space-time constitutional functions of polynucleotide chains; Informations by versifications

Part III Psychological Sensations and Perceptions

X Standpoint of the Psychological World

Definition of psychological perception; Empirically recognized postulates; Chromaticity plane; Vowel triangle; Projective-geometrical character of psychological sensations

XI Fundamental Structurology of Chromaticity

Barycentric coordinates; Equilateral barycentric representation; Horseshoe form of the spectrum locus; LAND's unorthodox colour synthesis

XII Contrast, Harmony and Anharmony between Quantized Psycho logical sensations

Finite fields; Triotron; tetratron; Cross ratio aspects of chromatic perceptions

XIII Vowel Harmony: Mathematical Study of Natural Language Structures – I

Triotronic harmony; Normal Vowel Harmony; Tetratonic inhibitions

XIV Quaternionian Phonoids and the Phonetical Tetrahedron

Fourth degree of freedom of information; Quaternionian character of the information propogations in the nerve system; Phonetical tetrahedron; Relation to the representation by the sensation function; Differentiation of quaternionian phonoids; Voice; Articulation

XV Spectrum Structure of Speech-sounds

Phonetical lattice; Contrast between palatal and labial articulations; A theoretical restriction on the domain of the sensation function

XVI Strata of Chromaticity Sensations

Development of cones and rods; Rod-cone break

XVII Consonant Vowels in Different languages: Mathematical Study of Natural Language Structures – 2

Consonants in the Chinese language; Germanic consonant shifts; The Centum and Satem shifts; Sanskrit vowel gradation; Phonetical structure of the ancient Japanese language

XVIII Natural Indo-European Grammar: Mathematical Study of Natural Language – 3

Design principle; Declension stem design; Declension ending design; Consonants in non-horizontal cycles; Dissection of the non-horizontal cycles; Dissection segments; Assignment to the Indo-European conjugation; Space- and time-constitution languages

XIX Galois Field Structures of Musical Scales

Definition of music; Algebra of musical scales; Structure of heptatonic; Harmony; Quantization by correlation of musical tones; Major and minor scales; Octave segmentation and leading note; Standpoint of the counterpoint; Oriental music and modern music; Practical materializations

**Appendix C Inaugral Address to the Einstein Centenary Symposium Nagpur, India, Feb 1980**

> *"Ladies and Gentlemen, President, Chairman, Secretary and Officers of the Organisation Committee and all participants. It is indeed a great honour and pleasure for me to be granted the privilege to address you inaugurating the four days for the Einstein Centenary Symposium now being open in this Central City of India, Nagpur.*
>
> *The Symposium may start with a recollection of how and why the late Albert Einstein is the most outstanding of all the great scientists this century has ever seen. One cannot sufficiently depict his incomparable standing, no matter how one may emphasize his great influences on contemporary sciences. Even by so doing many more things cannot help being missed. He combined in himself high ethical standards and marvellous achievements in scientific discovery and creation. All who gather today no doubt admit that Einstein is here with us. This great leader of science is indeed a friend of everybody he inspires. He is unique to have provided sources and stimuli for future developments of human cultural prospects. Not only science but philosophy also had to undergo revolutionary changes owing to the principle of relativity and revelation for quantum-physical premises opened by him.*
>
> *The world has celebrated these years many jubilees to felicitate Einstein. But India with her transcendency with philosophical surges from antiquity is most appropriate for reminding us of the profound epistemological implication of what Einstein taught us.*
>
> *After his great work, many tended to feel that there is hardly anything left to do. However, this is not the case. Einstein has left us tremendous prospect to encourage us. The framework of space-time, though not so much less than a hundred years old now, still provides the most powerful implement for penetrating the microscopic puzzles as well as permeating to the macroscopic world picture. Notwithstanding inadvertent debates that have not been lacking sometimes about mysterious matters, fields and spaces, unification is acquiring a more luminous aspect along the broad road originally paved by him. Which spirit other than his can one ultimately rely on in exploring the far region of the universe where black holes come about.*
>
> *Referring fundamental scientific problems that lie in all these it is assumed that Einstein did not show the final solution in his lifetime. What is most important is not specific, restricted formalisms but the structure of the truth that may have many faces. For exploring it, this Centenary Meeting has a most significant role to play from scientific, cultural standpoints working as a milestone implying a prodigious future prospect".*